\documentclass{IOS-Book-Article}

\usepackage{mathptmx}
\usepackage{graphicx}
\usepackage{amssymb,amsmath}
\usepackage{color}
\usepackage{hyperref}
\hypersetup{
  colorlinks   = true, 
  urlcolor     = blue, 
  linkcolor    = blue, 
  citecolor   = red 
}

\begin{document}
\begin{frontmatter}              

\title{A space-time parallel solver for the three-dimensional heat equation}
\runningtitle{IOS Press Style Sample}

\author[A,B]{\fnms{Robert} \snm{Speck}%
\thanks{E-mail: r.speck@fz-juelich.de}},
\author[B]{\fnms{Daniel} \snm{Ruprecht}},
\author[C]{\fnms{Matthew} \snm{Emmett}},
\author[D]{\fnms{Matthias} \snm{Bolten}}
and
\author[B]{\fnms{Rolf} \snm{Krause}}

\runningauthor{R.~Speck et al.}
\address[A]{J\"ulich Supercomputing Centre, Forschungszentrum J\"ulich, Germany}
\address[B]{Institute of Computational Science, Universit{\`a} della Svizzera italiana, Lugano, Switzerland}
\address[C]{Center for Computational Sciences and Engineering, Lawrence Berkeley National Laboratory, USA}
\address[D]{Department of Mathematics, Bergische Universit\"at Wuppertal, Germany}

%
%
%
%

\begin{abstract}
The paper presents a combination of the time-parallel ``parallel full approximation scheme in space and time'' (PFASST) with a parallel multigrid method (PMG) in space, resulting in a mesh-based
solver for the three-dimensional heat equation with a uniquely high degree of efficient concurrency. 
Parallel scaling tests are reported on the Cray XE6 machine ``Monte Rosa'' on up to $16\mathord{,}384$ cores and on the IBM Blue Gene/Q system ``JUQUEEN'' on up to $65\mathord{,}536$ cores.
The efficacy of the combined spatial- and temporal parallelization is shown by demonstrating that using PFASST in addition to PMG significantly extends the strong-scaling limit.
Implications of using spatial coarsening strategies in PFASST's multi-level hierarchy in large-scale parallel simulations are discussed.
\end{abstract}

\end{frontmatter}

\thispagestyle{empty}
\pagestyle{empty}

\section*{Introduction}
With the anticipated need for more than 100 million-way concurrency on emerging high-performance computing systems~\cite{IESR}, interest in additional ways to parallelize the solution of partial differential equations is rapidly increasing. 
For time-dependent problems, novel methods introducing concurrency along the temporal dimension are a promising and fascinating way to extend existing strong-scaling limits of space-parallel approaches.
A well-known time-parallelization method is Parareal, introduced in~\cite{LionsEtAl2001}, but other methods exist as well, e.g. the ``parallel implicit time algorithm'' (PITA)~\cite{FarhatEtAl2003} and the recent ``parallel full approximation scheme in space and time'' (PFASST)~\cite{MinionEtAl2008,Minion2010,EmmettMinion2012}.
PFASST is based on a combination of ``spectral deferred correction'' (SDC) methods~\cite{dutt_spectral_2000} and Parareal. 
It uses a hierarchy of space-time grids and includes a ``full approximation scheme'' (FAS) to enhance coarse-grid accuracy. 
The FAS correction allows different strategies for spatial coarsening to be used within the mesh hierarchy, leading to improved parallel efficiency.
An accuracy study of PFASST as well as serial SDC for a first-order particle-based discretization can be found in~\cite{SpeckEtAl2012_DDM}.
A successful demonstration of PFASST's efficacy in extreme-scale parallel particle simulations was presented in~\cite{SpeckEtAl2012}, where PFASST was combined with the parallel Barnes-Hut tree codes PEPC~\cite{winkel_massively_2012} to simulate a spherical vortex sheet.
However, to date, there are very few papers discussing performance of time-parallel methods in large- or extreme-scale parallel runs.\par
In the present paper, we introduce a massively space-time parallel solver for the three-dimensional heat equation based on the combination of PFASST with a mesh-based discretization of the Laplacian.
Here, space-time coarsening with PFASST is obtained using (i) a reduced set of temporal nodes, (ii) a hierarchy of spatial meshes, and (iii) different orders of the discretization of the Laplacian, following the strategies discussed in~\cite{SpeckEtAl2013_SISC}.
An implicit Euler method is used as sub-stepper within the SDC scheme, employing a parallel multigrid method (PMG)~\cite{chow2006} to solve the resulting linear problem.
The ability to efficiently solve problems of similar type to the heat equation constitutes a building block for tackling more complex problems. 
In the Navier-Stokes equations, for example, implicit-explicit methods are often applied where the diffusion term is treated implicitly and advection explicitly.
This leads to a structurally similar problem in the implicit step.
Thus, being able to efficiently solve the simpler problem studied here is in some sense a ``necessary condition'' to be able to solve more complicated setups, although more complex dynamics will likely result in less optimal convergence behavior of both PMG and PFASST.
We also show that running PFASST on top of a nearly saturated spatial parallelization has significant and interesting consequences for the coarsening strategies employed.

\section{Problem and Numerical Methods}\label{sec:setup}
The test problem considered here is the three-dimensional heat equation with a forcing term
\begin{align}
  \frac{\partial u(x,t)}{\partial t} &= \nu\Delta u(x,t) + f(x,t)\ \text{ in } \Omega = [0,1]^3,\\
  u(x,t) &= 0\ \text{ on } \partial\Omega
\end{align}
with $x=(x_1,x_2,x_3)^T\in[0,1]^3$, $t\in[0,T]$, the viscosity parameter $\nu = 0.1$ and $\Delta$ denoting the Laplacian. We choose 
\begin{align}
   f(x,t) = -\sin(\pi x_1)\sin(\pi x_2)\sin(\pi x_3)(\sin(t) - \nu\pi^2\cos(t))
\end{align}
for the source term, which leads to the analytic solution
\begin{align}
   u(x,t) = \sin(\pi x_1)\sin(\pi x_2)\sin(\pi x_3)\cos(t).
\end{align}
A method-of-lines approach is used, discretizing the problem in space first. 
Then, the time-parallel PFASST method, see~\S\ref{ssec:PFASST}, is applied to integrate the resulting system of initial value problems in time.
Here, PFASST uses a two-level hierarchy with $M=5$ collocation nodes per time step on the fine level and $M=3$ nodes per time step on the coarse level.
An implicit-explicit Euler method is used within PFASST to perform SDC sweeps, treating the diffusion term implicitly and the source term explicitly.
A parallel multigrid method, see~\S\ref{ssec:pmg}, is used to solve the linear systems of equations arising from the implicit part.
The problem is integrated until $T=6.0$ using steps of size $\Delta t = 0.1875$.
The Laplacian on the fine level is discretized using a fourth-order compact stencil with $255$ degrees-of-freedom per dimension, while on the coarse level, a classical second-order stencil is applied with $127$ degrees-of-freedom per dimension. 
Following the terminology introduced in~\cite{SpeckEtAl2013_SISC}, space-time coarsening is thus obtained by reduction of SDC nodes, by reduction of degrees-of-freedom in space and also by reduction of the spatial discretization order.
In all simulations, a threshold of $10^{-10}$ for the residual of the SDC iteration (time-serial) and PFASST iterations (time-parallel) is used, resulting in a relative maximum error of approx. $1.1\times 10^{-11}$.

\subsection{SDC and PFASST}\label{ssec:PFASST}
The time-parallel PFASST algorithm in its final form was introduced in~\cite{EmmettMinion2012} as a combination of spectral deferred correction methods~\cite{dutt_spectral_2000} with Parareal~\cite{LionsEtAl2001} plus an FAS correction to allow for efficient spatial coarsening along the level hierarchy.
It can also be understood as the time-parallel version of the recently introduced ``multi-level spectral deferred correction'' method (MLSDC, see~\cite{SpeckEtAl2013_SISC}).
Classical, single-level SDC methods employ a low-order method, typically an Euler scheme, to iteratively compute a collocation solution to an initial value problem.
Formally, using $k$ iterations or ``sweeps'' of a first-order method results in an order $\mathcal{O}(\Delta t^{k})$ SDC method, if the spectral collocation formula is accurate enough.\par
In SDC, a solution over a time step $[T^{n},T^{n+1}]$ with fixed length $\Delta t$ is computed based on the Picard formulation of the IVP 
\begin{equation}
  \label{eq:picard}
  u(T^{n+1}) = u(T^{n}) + \int_{T^{n}}^{T^{n+1}} f(u(s),s)\ \mathrm{d}s, \quad u \in \mathbf{R}^{N}
\end{equation}
using a collocation formula with sub-steps $T^{n} = t_{0} < t_{1} < \ldots < t_{M} = T^{n+1}$ and weights $q_{m,j}$ to approximate~\eqref{eq:picard} by 
\begin{equation}
  \label{eq:colloc}
  U_{m} = U_{0} + \Delta t \sum_{i=0}^{M} q_{i,m} f(U_{i},t_{i})
\end{equation}
with $U_{m} \approx u(t_{m})$.
Typically, Gauss-Lobatto nodes are used so that $T^n=t_0$ as well as $T^{n+1}=t_M$ are part of the set of collocation nodes.
Equation~\eqref{eq:colloc} corresponds to a $N \times (M+1)$ nonlinear system of equations for the unknowns $U_{m}$. 
Directly applying some suitable solver would yield a collocation method, a subclass of implicit Runge-Kutta methods, see e.g.~\cite{HairerWanner2006}.
In SDC, instead of directly solving the full system, a low-order method is used to generate an iterative solution that converges (usually exponentially fast) to the collocation solution.
For an IMEX-Euler based sub-stepping scheme (which uses an implicit Euler stepper for the implicit/stiff part $f^{\rm I}$ of $f$ and an explicit Euler stepper for the explicit/non-stiff part $f^{\rm E}$ of $f$), one sweep corresponds to computing
\begin{align}
  U^{k+1}_{m+1} = U^{k+1}_{m} &+ \Delta t_{m} \left[ f^{\rm E}( U^{k+1}_{m}, t_{m} ) - f^{\rm E}(U^{k}_{m}, t_{m}) \right] \nonumber \\
			      &+ \Delta t_{m} \left[ f^{\rm I}( U^{k+1}_{m+1}, t_{m+1}) - f^{\rm I}(U^{k}_{m+1}, t_{m+1}) \right] \nonumber \\
			      &+  \Delta t S^{k}_{m}
\end{align}
for $m=0,\ldots,M-1$ where $k$ denotes the SDC iteration number.
Here, $\Delta t_{m} = t_{m+1} - t_{m}$ denotes the distance between two collocation points while
\begin{equation}
  \Delta t S^{k}_{m} \approx \int_{t_{m}}^{t_{m+1}} f(U^{k}(s),s)\ \mathrm{d}s
\end{equation}
using a suitable spectral quadrature rule.
If the iteration converges, that is if $U^{k+1}_{m} - U^{k}_{m} \to 0$ as $k \to \infty$, the solution at the last collocation point $t_{M} = T^{n+1}$ becomes
\begin{equation}
  U^{k+1}_{M} = U^{k+1}_{0} + \Delta t \sum_{m=0}^{M} S^{k}_{m} \approx u(T^{n}) + \int_{T^{n}}^{T^{n+1}} f(u(s),s)\ \mathrm{d}s = u(T^{n+1})
\end{equation}
and thus yields an approximate solution at $T^{n+1}$ identical to what would be obtained from directly computing the collocation solution.\par
In multi-level SDC, instead of computing sweeps on one level, sweeps are performed on a hierarchy of levels where higher (i.e. coarser) levels use fewer collocation points and usually also a coarsened spatial representation of the problem, see~\cite{SpeckEtAl2013_SISC} for details.
In PFASST, iterations of MLSDC are performed concurrently on multiple time slices.
Updated values are passed forward in time as soon as a sweep on a particular level has finished, which subsequently serve as new initial values for the following processor/time step.
However, blocking communication is required on the coarsest level only~\cite{EmmettMinion2012_DDM}, so that PFASST requires minimal synchronicity between the multiple MLSDC iterations.
The reader is referred to~\cite{EmmettMinion2012,SpeckEtAl2012,SpeckEtAl2013_SISC} for details.
The speedup to be expected from two-level PFASST on $P_{T}$ many processors on an identical number of time steps compared to a serial SDC run can be modeled as
\begin{equation}
  \label{eq:pfasst_speedup}
  s(P_{T}) = \frac{ K_{S} P_{T} }{ P_{T} \alpha + K_{P} \left( 1 + \alpha + \beta \right)}
\end{equation}
see~\cite{Minion2010,SpeckEtAl2012}. 
Here, $\alpha$ is the ratio of the runtime required for one coarse sweep to the runtime for one fine sweep, $K_{S}$ is the number of sweeps required by the
serial SDC method, and $K_{P}$ is the number of PFASST iterations required to converge to the same accuracy as the serial SDC scheme. 
Finally, $\beta$ denotes overhead e.g. from interpolation, restriction or communication in PFASST.
Note that minimizing the ratio $\alpha$ without increasing the number of required iterations $K_{P}$ too much is key to achieving good parallel speedup with PFASST.
Different strategies are identified and discussed in~\cite{SpeckEtAl2013_SISC} in the context of MLSDC, which can be readily applied for PFASST as well.
For suitable setups, the number of iterations is reduced by applying MLSDC instead of SDC, since sweeps on the coarse (and computationally cheap) level are partially able to replace fine sweeps.
This effect is somewhat lessened in the case of PFASST.
Here, depending on the problem, the concurrent iteration over multiple time steps naturally leads to an increased number of iterations required for convergence, since early iterations on later time steps start with rather inaccurate initial values.

\subsection{Parallel Multigrid in Space}\label{ssec:pmg}
For solving the linear systems arising in each step of the implicit-explicit SDC scheme and for inverting the weighting matrix to compute the FAS correction, a parallel multigrid method (PMG) is used.
Multigrid methods are optimal or near optimal for a wide range of problems. 
A general introduction to parallel multigrid methods can be found for example in~\cite{chow2006}.
PMG distributes the unknowns among the processors by splitting the domain.
On coarser grids unknowns stay on the processor they have been assigned to on the finest grid, resulting in idle processors. The coarsest grid consists of one unknown.
For more details on PMG and its inclusion into MLSDC and PFASST, particularly in combination with compact finite difference stencils, the reader is referred to~\cite{SpeckEtAl2013_SISC}.

\section{Performance Results}

The setup described in \S\ref{sec:setup} is run on two different machines: The IBM Blue Gene/Q ``JUQUEEN'' at J\"ulich Supercomputing Centre and the Cray XE6 ``Monte Rosa'' at the Swiss National Supercomputing Centre.
JUQUEEN features $28\mathord{,}672$ nodes, each equipped with an 16-core IBM PowerPC A2 1.6GHz processor, yielding a total of $458\mathord{,}752$ cores and a peak performance of $5.9$ PFlops (Linpack: 5.0 PFlops).
The nodes, each with $16$ GB main memory, are connected by a 5D torus interconnect.
Monte Rosa consists of 1496 nodes, each with two 16-core AMD Interlagos Opteron CPUs with $2.1$ GHz for a total of $47\mathord{,}872$ cores and a peak performance of $402$ TFlops (Linpack: 319 TFlops).
Each node has $32$ GB of memory and nodes are connected by a Gemini 3D interconnect.
On both machines, one small and one large run of PMG+PFASST is conducted where the small runs use a number of cores for PMG that corresponds to the point where PMG's scaling saturates for the coarse level problem and the large runs use a number of cores that correspond to the saturation level of the fine level problem.
For comparison, two serial SDC runs using only PMG for spatial parallelization are also performed, one solving the fine level problem the other solving the coarse level problem.

\subsection{Space-time parallel speedup on the IBM Blue Gene/Q}\label{ssec:scaling_Q}
\begin{figure}[t]
   \centering
   \includegraphics[width=\textwidth]{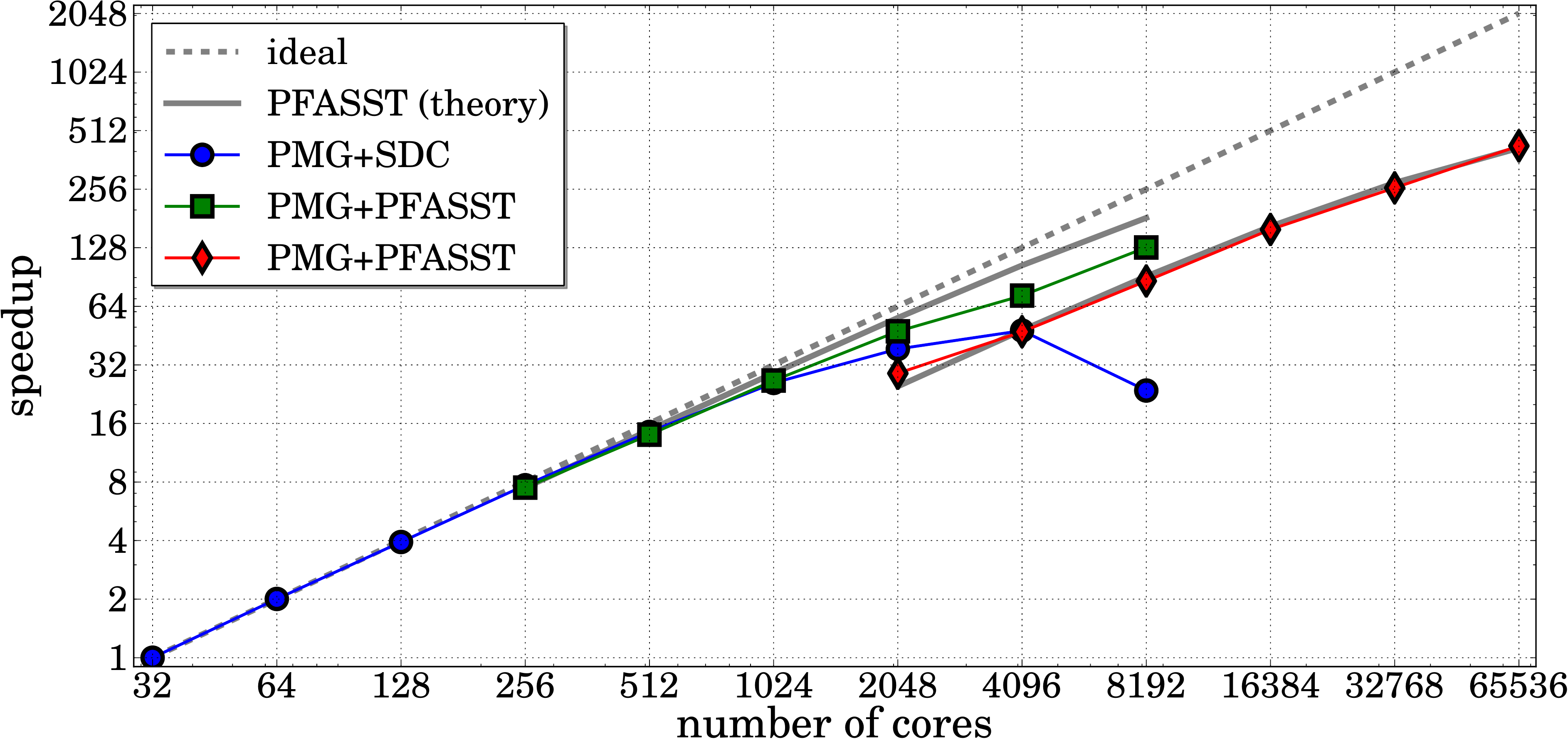}
  \caption{Speedup of time-serial PMG+SDC (blue, circle marker) as well as PMG+PFASST with PMG using $256$ cores (green, square marker) or $2\mathord{,}048$ cores (red, diamond marker) on the IBM Blue Gene/Q JUQUEEN.
The dashed grey line indicates ideal speedup, the solid gray lines the speedup from PFASST according to the theoretical estimate~\eqref{eq:pfasst_speedup}.}\label{fig:scaling_bgq}
\end{figure}
Figure~\ref{fig:scaling_bgq} shows the scaling of the code on the IBM Blue Gene/Q system. 
The blue line with circle markers indicates the speedup of PMG alone, that is without temporal parallelization, solving the fine level problem with plain single-level SDC time stepping.
The dashed gray line marks ideal speedup.
For the chosen problem, PMG scales near perfectly to $1\mathord{,}024$ cores and continues to scale with less than ideal efficiency to $4\mathord{,}096$ cores.
Beyond $4\mathord{,}096$, using more cores reduces speedup.
For comparison, another serial PMG run solving the coarse level problem is performed (not shown here), indicating that for the smaller coarse problem PMG scales up to $256$ cores.
Note that this strong-scaling capability is by no means optimal for linear multigrid solvers.
In the present version of the code, being a more ``black-box``-type coupling of PMG to the PFASST library, optimization is still necessary to increase the performance of the spatial solver.
Preliminary tests not documented here suggest that the sub-optimal scaling of PMG within SDC and PFASST stems from a load balancing issue caused by multiple calls of the solver in rapid succession.
Nevertheless, the results obtained here provide a compelling illustration of the potential as well as valuable insight into the limitations of space-time parallel solvers, see\S\ref{ssec:efficiency} for more discussion.\par
Two simulations are performed with PMG + PFASST, one using 256 cores for PMG times $2$, $4$, $8$, $16$ and $32$ cores for PFASST (green line, square markers),
the other one using $2\mathord{,}048$ cores for PMG and also $2$, $4$, $8$, $16$ and $32$ cores for PFASST (red line, diamond markers).
Speedup is measured against the serial, single-level execution of the code with $32$ cores in space, i.e.~against the SDC results with $32$ cores for PMG.
The solid gray line marks the speedup according the the theoretical estimate for PFASST.
The parameter $\alpha$ in~\eqref{eq:pfasst_speedup} is computed numerically from the two serial PMG runs, yielding
\begin{equation}
  \alpha_{\rm small} \approx 0.083 \quad \text{and} \quad \alpha_{\rm large} \approx 0.36
\end{equation}
on JUQUEEN.
Note that $\alpha$ is about a factor $4.3$ larger in the large run than in the small run and thus the speedup to be gained from PFASST is noticeably smaller there.
The reason is that if PMG uses sufficiently many cores to reach its scaling limit, the employed reduction of degrees-of-freedoms is significantly less efficient, see the discussion in \S\ref{ssec:efficiency}.\par
The runs using $2\mathord{,}048$ cores for PMG correspond to the saturation point of the fine level problem while runs using $256$ cores correspond to the saturation point of the coarse level problem.
Note that the first data point of the red line at $2\mathord{,}048$ cores indicates a run with a single time-rank, for which PFASST is identical to MLSDC.
The drop in speedup compared to single-level SDC (blue line) stems from the fact that MLSDC is somewhat more expensive than SDC for the setup studied here.
Similarly, the first data point of the green line also corresponds to MLSDC, but here runtimes of both single-level SDC and MLSDC are essentially identical.
For both cases, as the number of time-ranks increases, strong scaling of the space-time parallel PMG+PFASST combination goes significantly further than that of PMG alone.
The small run using fewer cores for PMG, however, consistently provides better speedup than the large run for an identical number of cores.\par
Table~\ref{tab:eff_bgq} summarizes the parallel efficiency of the temporal parallelization of both runs, that is of the additional speedup provided by PFASST.
The less then unity speedup indicated for a single time-rank for both cases quantifies the overhead of using MLSDC over SDC discussed above.
Clearly, the small run show significantly better efficiency, again because of a better optimized temporal parallelization with a much smaller $\alpha$.
\begin{table}[t!]
\begin{tabular}{ccc||ccc}
\multicolumn{6}{c}{\bf IBM Blue Gene/Q JUQUEEN} \\
\multicolumn{3}{c}{Small run with 256 cores PMG (green)}    & \multicolumn{3}{c}{Large run with 2048 cores PMG (red)}  \\ 
Time-ranks & Speedup & Efficiency & Time-ranks & Speedup & Efficiency \\ \hline
1                & 0.97        &      ---       &     1             & 0.75       &  ---                \\           
2                & 1.82        & 91.0\%    &     2             & 1.23       &  61.5\%               \\           
4                & 3.45        & 86.2\%    &     4             & 2.24       &  56.0\%               \\           
8                & 6.18        & 77.3\%    &     8             & 4.11       &  51.4\%               \\           
16              & 9.43        & 58.9\%    &     16           & 6.73       &  42.1\%               \\           
32              & 16.68      & 52.1\%    &    32           & 11.06     &  34.6\%               \\           
\end{tabular}\vspace{0.5em}
\caption{Efficiency of the speedup provided by PFASST. Here, the reference is the serial SDC run with $256$ cores (small run, $129.04$ sec for a single SDC time step) and $2\mathord{,}048$ cores (large run, $25.73$ sec for a single SDC time step) for PMG. Using PFASST, the small run was performed in $247.61$ sec on $32$ time steps simultaneously, while the large run was done in $74.44$ sec.}
\label{tab:eff_bgq}
\end{table}

\subsection{Space-time parallel speedup on the Cray XE6}\label{ssec:scaling_R}
\begin{figure}[t]
   \centering
   \includegraphics[width=\textwidth]{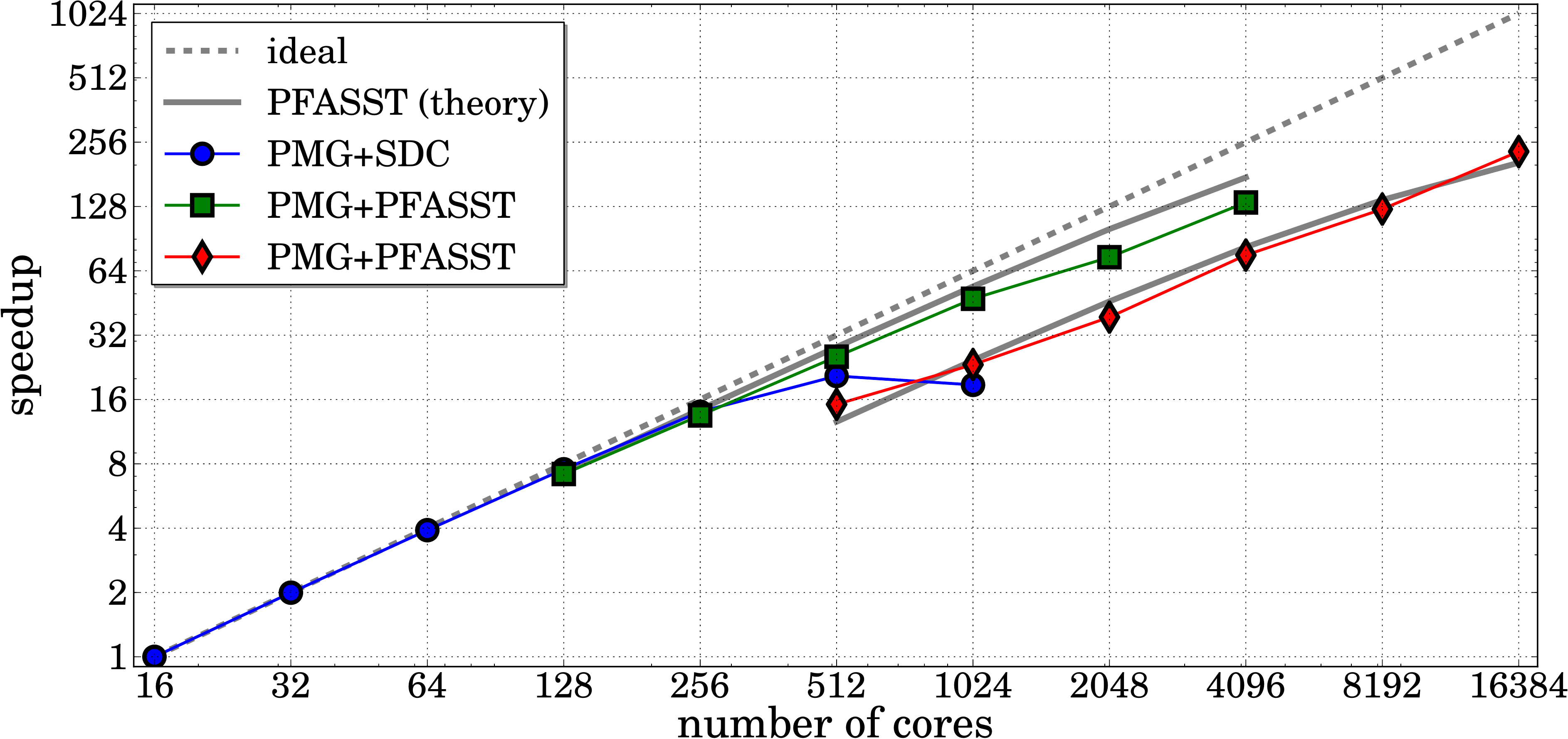}
  \caption{Speedup of PMG+SDC (blue, circle marker) as well as PMG+PFASST with PMG using $128$ cores (green, square marker) or $512$ cores (red, diamond marker) on the Cray XE6 Monte Rosa.
The dashed gray line indicates ideal speedup, the solid gray lines the speedup from PFASST according to the theoretical estimate~\eqref{eq:pfasst_speedup}.}
\label{fig:scaling_rosa}
\end{figure}

As on JUQUEEN, results from three runs performed on the Cray XE6 are shown in Figure~\ref{fig:scaling_rosa}: One with SDC+PMG and no temporal parallelization,
one small run of PMG+PFASST using as many cores for PMG as are necessary to saturate it for the coarse problem (green line, square markers) and
one large run, using the number of cores corresponding to the saturation limit of PMG for the fine problem (red line, diamond markers).
Furthermore, a serial run of PMG with the coarse level spatial discretization is done for comparison (not shown here).
The dashed line indicates ideal speedup.\par
On the Cray, PMG scales almost perfectly to $256$ and saturates at $512$ cores for the fine level problem.
Beyond 512 cores, runtimes go up as in \S\ref{ssec:scaling_Q}.
The scaling limit for PMG when solving the coarse problem with $127^{3}$ mesh points (again not shown) is reached at $128$ cores.
The solid line indicates the speedup to be expected from PFASST according to estimate~\eqref{eq:pfasst_speedup}.
As before, the parameter $\alpha$ is computed numerically from the serial PMG runs, yielding
\begin{equation}
  \alpha_{\rm small} \approx 0.086 \quad \text{and} \quad \alpha_{\rm large} \approx 0.4
\end{equation}
on Monte Rosa.
Again, the parameter is much higher in the large runs (by a factor of $4.7$), indicating substantially less efficient parallelization by PFASST.
For both the small and large run, PMG+PFASST provides significant additional strong scaling.
Both runs also match the theoretical predictions, but again the smaller run using fewer cores for PMG shows significantly better speedups for identical numbers of cores.
Table~\ref{tab:eff_rosa} shows the efficiency achieved on the Cray and as before, the smaller run turns out to be significantly more efficient.
\begin{table}[ht]
\begin{tabular}{ccc||ccc}
\multicolumn{6}{c}{\bf Cray XE6 Monte Rosa} \\
\multicolumn{3}{c}{Small run with 128 cores PMG (green)}    & \multicolumn{3}{c}{Large run with 512 cores PMG (red)}  \\ 
Time-ranks & Speedup & Efficiency & Time-ranks & Speedup & Efficiency \\ \hline
1                & 0.95        &      ---       &     1             & 0.73       &  ---                \\           
2                & 1.79        & 89.5\%    &     2             & 1.13       &  56.6\%               \\           
4                & 3.36        & 84.0\%    &     4             & 1.89       &  47.2\%               \\           
8                & 6.28        & 78.5\%    &     8             & 3.68       &  46.0\%               \\           
16              & 9.82        & 61.4\%    &     16           & 6.00       &  37.5\%               \\           
32              & 17.72      & 55.4\%     &    32           & 11.22     &  35.1\%               \\           
\end{tabular}\vspace{0.5em}
\caption{Efficiency of the speedup provided by PFASST. Here, the reference is the serial SDC run with $128$ cores (small run, $73.42$ sec for a full SDC time step) and $512$ cores (large run, $26.88$ sec for a full SDC time step) for PMG. With PFASST, $32$ parallel time steps took $132.09$ sec with the small run and $76.64$ sec with the large run.}
\label{tab:eff_rosa}
\end{table}

\subsection{Efficiency of spatial coarsening at saturation point}\label{ssec:efficiency}
On both systems, the IBM Blue Gene/Q and the Cray XE6, the space-time parallel method achieves larger speedup and better efficiency if the number of cores for the spatial parallelization is not set so high as to fully saturate PMG for the fine level problem.
This is because the speedup from PFASST depends critically on the parameter $\alpha$ in~\eqref{eq:pfasst_speedup}, which denotes the runtime ratio of coarse to fine sweeps.
If the efficiency of the spatial parallelization is sufficiently low so that runtimes are strongly dominated by communication time, reducing the degrees-of-freedom on the coarse level does not yield a significant reduction in runtime for the coarse sweeps.
Although runtimes may decrease, the dominating communication time stays approximately the same, and thus the difference in runtime for sweeps on the coarse and fine level remains small.
This is demonstrated by observing that $\alpha$, which ideally should be close to zero, was more than a factor of four times larger for the big runs than for the small runs on both machines:
if a smaller number of cores is used in space, the smaller problem on the coarse level achieves significantly reduced runtimes compared to the fine level.
As shown by the results above, this strategy results in much better efficiency of the temporal parallelization and better overall speedups of the combined space-time parallel method.
In regimes where the efficiency of PMG is far from optimal, the additional runtime associated with using fewer cores for PMG is more than compensated for by the improved efficiency of using PFASST+PMG.
This illustrates the need to consider \emph{all} factors determining the runtimes in estimates for the speedup of time parallel methods in large-scale parallel simulations.
Accounting only for pure computing times will lead to non-optimal strategies.\par
It is interesting to note that a similar problem is encountered in parallel multigrid methods when the number of degrees-of-freedoms on the coarse levels in the hierarchy becomes too small, see~\cite{HuelsemannEtAl2006}.
Adopting strategies developed for space multigrid methods to cope with this problem, e.g. ``coarse grid agglomeration'', could be a successful approach to optimize combined space-time parallelizations.\par
These insights also emphasize that using time-parallel methods on top of spatial parallelization to extend strong scaling is not straightforward in the sense of simply adding the former to the latter:
in large-scale parallel runs, temporal and spatial parallelization strategies interact and the properties of one can directly affect the performance of the other.

\section{Conclusions}
This paper presents a study of the performance of a combination of the time-parallel ``parallel full approximation scheme in space and time'' (PFASST) with a parallel multigrid solver in space, resulting in a space-time parallel solver for the three-dimensional heat equation.
Runtimes are reported on up to $16\mathord{,}384$ cores on a Cray XE6 as well as on up to $65\mathord{,}536$ cores on an IBM Blue Gene/Q system.
Speedups on the latter are somewhat better, but the differences are small and in general theoretical estimates for PFASST are found to be quite accurate.
It is shown that PFASST provides significant additional strong scaling beyond the saturation of the spatial solver.\par
To improve its parallel efficiency, PFASST employs a hierarchy of space-time meshes (two in the present study): decreasing the relative runtime on the coarse level compared to the fine level results in better efficiency.
To reduce the runtime for sweeps on the coarse level, the order and number of mesh points of the spatial discretization were reduced in conjunction with temporal coarsening.
The results obtained here suggest that the efficiency of using this strategy to optimize the scaling of PFASST is not straightforward for massively parallel simulations.
In particular, when the spatial parallelization is strongly dominated by communication time, reducing the degrees-of-freedom on the coarse level does not significantly reduce the runtime of coarse sweeps, which in turn limits the efficiency of time-parallelization.  In contrast, using fewer cores in space yields a more optimal coarse-to-fine runtime ratio and significantly improves the efficiency of PFASST.
The improvement brought by PFASST offsets the increased runtimes resulting from using fewer cores for the spatial solver and, overall, leads to better speedup of the space-time parallel approach.
This demonstrates that when parallel-in-time methods are considered in the context of massively parallel simulations as a means of extending the strong scaling of space-parallelization, the interplay between space- and time-parallelization must be taken into account.

\paragraph{Acknowledgments.} 
Computing time on Monta Rosa at the Swiss National Supercomputing Centre (CSCS) was provided by project D16 while
computing time on JUQUEEN at J\"ulich Supercomputing Centre (JSC) was provided by project HWU12.
D.~Ruprecht and M.~Emmett also thankfully acknowledge support by the Swiss National Science Foundation (SNSF) through Grant 147597.
We also thank M.~Winkel from JSC for compiling a modern version of GCC on the IBM Blue Gene/Q.
\bibliographystyle{unsrt}
\bibliography{parco2013}

\end{document}